\documentclass[final,authoryear]{autart}

\usepackage{amsfonts}
\usepackage{enumerate} \usepackage{verbatim} \usepackage{amssymb}
\usepackage{setspace} \usepackage{epsfig} \usepackage{url} \usepackage{graphicx} \usepackage{subfigure}
\usepackage{natbib}

\usepackage{amsmath} \usepackage[amsmath,amsthm,thmmarks]{ntheorem}

\def\urltilde{\kern -.15em\lower .7ex\hbox{\~{}}\kern .04em}
\def\urldot{\kern -.10em.\kern -.10em}

\def\urlhttp{http\kern -.10em\lower -.1ex\hbox{:}\kern -.12em\lower 0ex\hbox{/}\kern -.18em\lower 0ex\hbox{/}}

 \newtheorem{theorem}{Theorem}
 \newtheorem{corollary}{Corollary}
  
  { \theoremsymbol{\rule{1ex}{1ex}}  \newtheorem{example}{Example}  }
 \newtheorem{problem}{Problem}
 \newtheorem{proposition}[theorem]{Proposition}
 \newtheorem{remark} {Remark}

\newcommand {\R}{\mathbb R}

 \newcommand {\SA}{\mathcal S}

 \newcommand{\sname}{}
\newcommand{\slabel}[1]{\debug{\fbox{\tiny \sname #1}}\label{\sname #1}}
\newcommand{\debug}[1]{}              % final version
\newcommand{\FB}{\begin{figure}[t]\centering} \newcommand{\LFB}{\begin{figure*}[t]\centering}
\newcommand{\FE}[2]{\caption{#2 \debug{\fbox{\sname #1}}} \slabel{#1} \end{figure}}
\newcommand{\LFE}[2]{\caption{#2 \debug{\fbox{\sname #1}}} \slabel{#1} \end{figure*}}
\newcommand{\tB}{\begin{table}[hbtp]\centering} \newcommand{\tE}[2]{\caption{#2 \debug{\fbox{\sname
#1}}}\slabel{#1} \end{table}} 

\newcommand{\be}{\begin{equation}}
\newcommand{\ee}{\end{equation}}
\begin{document}
 \begin{frontmatter}

\title{On  Boolean Control Networks with Maximal Topological~Entropy\thanksref{footnoteinfo}}
%%%
\thanks[footnoteinfo]{This paper was not presented at any IFAC
meeting.
The research of MM is supported in part by a research grant from the Israel Science Foundation~(ISF).
Corresponding author: Prof. Michael Margaliot,
  Tel: +972 3 640 7768;
Email: \texttt{michaelm@eng.tau.ac.il}}
\author[First]{Dmitriy Laschov}
\author[Second]{Michael Margaliot}

\address[First]{School  of Elec. Eng.-Systems, Tel Aviv University,  Israel.}
\address[Second]{School  of Elec. Eng.-Systems, Tel Aviv University,  Israel.}

   \maketitle

\begin{abstract}
%%%%%%%%%%%%%%%%%
 Boolean control networks~(BCNs) are discrete-time
 dynamical systems with    Boolean state-variables and   inputs
 that are interconnected via
 Boolean functions.
 BCNs are recently attracting considerable interest  as
computational models for genetic and cellular
 networks with exogenous inputs.

 The topological entropy  of a BCN  with~$m$ inputs is a
nonnegative real number in the  interval~$[0,m \log 2]$.
Roughly speaking, a larger topological entropy
means that asymptotically the control   is    ``more powerful''.
We derive a necessary and sufficient condition
for a BCN to have the  maximal possible topological entropy.
Our condition is stated in the framework of Cheng's algebraic state-space representation
of BCNs. This means that verifying this condition incurs an exponential time-complexity.
We also show that the problem of determining whether a BCN with~$n$ state variables and~$m=n$ inputs
has a maximum   topological entropy is NP-hard,
suggesting that this problem cannot be solved in general
using a polynomial-time algorithm.
%%%%%%%%%%%
\end{abstract}

\end{frontmatter}

\section{Introduction}
%%%%%
%%%%%%%

 Boolean   networks~(BNs) are useful modeling tools for   dynamical systems whose  state-variables
    can attain  two possible values. Examples range from
   artificial neural networks with ON/OFF type
  neurons   (see, e.g.~\cite{hass95}),
to models for the
   emergence of social consensus between  simple  agents that can either agree or
  disagree with a certain opinion
(see, e.g.~\cite{Green}).

There is a growing interest in modeling
biological systems using BNs and, in particular,
 genetic regulation networks, where each gene can be
  either
  expressed~(ON) or not
  expressed~(OFF)~(\cite{chaos,kauffman2003,li_yeast}).
Although being highly abstract,
BNs seem to capture the real behavior
 of gene-regulatory processes well~(\cite{bornholdt,10.1109/TCBB.2011.62})

\cite{kauff69} has
  studied the  order and stability of
  large, randomly constructed nets of such binary genes.
 He also  related the behavior of these random nets to various
cellular control processes, including  cell differentiation,
by associating
 every  possible cell type
  with  a  stable attractor of the~BN.
%%%
This work has stimulated   the analysis
  of  large-scale
BNs
  using  tools from  the theory of  complex systems and statistical physics
(see, e.g.~\cite{albert,aldana,drossel,kauff93}).

% There are several other motivations for using
%BNs in systems biology~\cite{huang_2002}, including the fact that many metabolic and
%genetic networks demonstrate
%some form of bi-stability. An
%important  example  are  epigenetic
% switches (see, e.g.~\cite{ptashne}).

BNs have also been used to model various cellular processes
including the  complex cellular signaling network controlling stomatal closure
in plants~(\cite{assmann}),
 the molecular pathway between two neurotransmitter
systems, the dopamine and glutamate receptors~(\cite{gupta}), carcinogenesis,  and the effects   of
therapeutic intervention~(\cite{szallasi}).
%(see also~\cite{kauff71}).

%Modeling using BNs requires   coarse-grained qualitative information
%(e.g., the effect of one gene on another gene     is  either  activating or inhibiting).
%This is in sharp contrast  to  other models,
% for example,
% those based on   differential equations, that
% entail  specifying    numerous parameter values
%(e.g., rate constants). Several studies have  developed algorithms for
%deriving computational models in the form of BNs from
%  biological data  (see e.g.~\cite{akutsu_identify,multi_bin,layek2011} and the references
%therein).
%
%These studies suggest   that   BNs  provide a
%a  powerful
% tool for modeling large-scale biological networks~\cite{born08}.
%Computational models based on BNs are able to reproduce the
%main characteristics of the biological network dynamics:
%attractors of the BN correspond to stationary biological states;
%large attraction basins indicate   robustness of the biological state; and so on.

%For a general  exposition on various
% approaches for   modeling gene regulation networks, see~\cite{GRN_primer,shamir_nat}.

%Modeling a biological system
%involves  considerable uncertainty,  due to   perturbations that
%affect the biological system  and      inaccuracies of  the  measuring
%equipment.
%Incorporating  this  uncertainty in the modeling stage leads to
%\emph{Probabilistic Boolean Networks}~(PBNs)~\cite{shmulevich}.
%A PBN is a collection of (deterministic) BNs
%combined with a probabilistic switching rule determining which
%network  is active at each time instant.

BNs with (Boolean) inputs are referred to as Boolean control networks~(BCNs).
BCNs have been used to model biological systems with exogenous inputs.
For example,
 \cite{faure06} (see also~\cite{faure_2009}) have developed a BCN model for
the core network regulating the mammalian cell cycle.
Here the nine state-variables represent
the activity/inactivity of nine different proteins: Rb, E2F, CycE, CycA, p27, Cdc20, Cdh1, UbcH10, and CycB,
and
the single Boolean
input represents the activity/inactivity of CycD in the cell.
%%%

\cite{cheng_book}  have developed an
\emph{algebraic  state-space representation}~(ASSR)
 of      BCNs (and, in particular, of BNs).
 This representation has
  proved
    useful for studying
  control-theoretic questions, as they reduce a BCN to
a positive linear switched system whose input, state and output variables are canonical vectors.
Topics that have been analyzed using the ASSR include
optimal control~(\cite{Optimal_Control_cheng,dima_mul,dima_minimum_time}), controllability and observability~(\cite{dima_cont,Fangfei,cont_observ_bool_cont_net,EF_MEV_BCN_obs2012}),
 identification~(\cite{Cheng_Identification}),
disturbance decoupling~(\cite{distr_decoupling_BCN}), and more.

The ASSR of a BN with~$n$ state-variables and~$m$ inputs
  includes a~$2^n \times 2^{n+m}$
matrix. Thus, any algorithm based on the ASSR has an exponential time complexity. A natural question
is  whether better algorithms exist.
 \cite{zhao_remark} has
 shown   that
determining whether a  BN has a fixed point   is
  NP-complete.
 \cite{bcn_complexity} have   shown that several control problems for BCNs are  NP-hard.
 \cite{dima_obs} have
 shown  that the observability problem for BCNs is also NP-hard.
 Thus, unless $P=NP$,   these analysis problems for BCNs cannot be solved
in  polynomial time.
%For a general survey on the computational complexity of problems in systems and control theory,
%see~\cite{blondel_2000}.

\cite{hochma_sym}
noted the connection between~BCNs and \emph{symbolic dynamics}~(SD).
The main object of study in SD is
\emph{shift spaces} (\cite{Lind_1995}).
The set of all possible trajectories of a  BCN is  a shift space, so
  many results and analysis tools from SD are immediately applicable to BCNs.
  In particular, \cite{hochma_sym} noted
  that an important
   notion from SD called topological entropy can be defined for BCNs,
  and computed using the Perron root of a certain non-negative matrix that appears
  in the ASSR of a BCN.
  The topological entropy of a BCN with~$n$ state-variables and~$m$ inputs
  (we always assume that~$m\leq n$) is a number in the range~$[0,m\log 2]$ that
    indicates how ``rich'' the control is.

   In this paper, we derive a necessary and sufficient condition for a BCN to have
   a maximal topological entropy. This condition is stated in terms of the ASSR.
   We also show that for a BCN with~$n$ state variables and~$m=n$ inputs
   the problem of determining whether the  BCN has maximal  topological entropy is NP-hard.
This implies
 that unless~$P=NP$, there does not exist an algorithm with polynomial time complexity that solves this problem.

The remainder of this note  is organized as follows. Section~\ref{sec:preli} reviews BNs, BCNs,
  and some definitions and tools from SD.
 Section~\ref{sec:main}   includes  our main results.
  Section~\ref{sec:concl}   concludes and describes some possible directions for further research.

\section{Preliminaries} \label{sec:preli}
%%%%%%%%%%%%%%%%%%%%%%%%%%%%%%%%%%%%%%%%%%%%%%
%%%%%%%%%%%%%%%%%%%%%%%%%%%%%%%%%%%%%%%%%%%%%

We begin by reviewing   BCNs and their  ASSRs.
 Let~$ \SA:=\{0,1 \}$.
A BCN  is a discrete-time logical dynamical
  system
\begin{align}\label{eq:BCN}
%%%%%55
                    X_1(k+1)&=f_1(X_1(k), \dots, X_n(k),U_1(k),\dots,U_m(k)  ),\nonumber  \\
                             & \vdots \\
                    X_n(k+1)&=f_n(X_1(k), \dots, X_n(k) , U_1(k),\dots,U_m(k) ), \nonumber
 %%%%%
\end{align} where~$X_i ,U_i  \in \SA$,
 and each~$f_i$   is a Boolean  function, i.e. $f_i: \SA^{n+m} \to \SA $.
 It is useful to write this   in vector form as
 \be \label{eq:bcnio}
            X(k+1)=f(X(k),U(k)).
 \ee
%\begin{remark}\label{rem:FSM}
%Note that the BCN~\eqref{eq:bcnio} is a~\emph{finite-state machine} (FSM). And %since output~$y(t)$
%depends only on state~$x(t)$, BCN is~\emph{Moore machine}.
%\end{remark}
%%
A BN is a
  BCN without inputs, i.e.
  \be \label{eq:BNbinary}
                                X(k+1)=f(X(k)).
\ee
%
% It is worth noting that a BCN with~$m$ inputs may be interpreted as
%a \emph{switched system} switching between~$2^m$ possible  subsystems,
%where each subsystem is a BN, and with
%  the value of the control determining which subsystem is active at
%every time step. The next example demonstrates this.
%%%%%%%%%%%%%%%%%%%%%%%%%%%%%%%%%%%%%%%%%%%%%%%%%%%%%%%%%%%%%
%\begin{example}\label{exa:simple}
%%%%
%Consider the two-state, two-input  BCN
%%%%%
%\begin{align}   \label{eq:dynsex}
%%%%%%
%                    X_1(k+1)&= X_1(k)\vee [X_2(k)\wedge  U_1(k) ],\\
%                    X_2(k+1)&= X_2(k) \wedge            U_2(k).\nonumber
%     %%%
%%%%%%
%\end{align}
%%%%
%The control variables may attain one of
%four possible values:~$(U_1(k),U_2(k))\in\{ 11,10,01,00 \}$. With each of these
%values we can associate a corresponding subsystem which is a BN.
%For example, when   $U_1(k)=U_2(k)=1$   the corresponding
%dynamical equations are
%%
%\begin{align*}
%%%%%%
%                    X_1(k+1)& =X_1(k)\vee X_2(k)  ,\\
%                    X_2(k+1)&=X_2(k).
%%%%%%
%\end{align*}
%\end{example}
%
%
%In the ASSR reviewed below, each subsystem becomes
%a positive linear system, so a BCN becomes  a discrete-time
%positive linear switched system~(PLSS).
%For more on PLSSs, see e.g.~\cite{lior,mason-shorten-cdc04, valcher_sicon_2011,lior_SIAM}
%and the references therein.
%
%
%
%\subsection{Algebraic representation of BCNs}
%
\cite{cheng_book} have
    developed an
algebraic state-space representation of BCNs using the  semi-tensor product of matrices.
This topic has  been described in many publications, so we review it briefly.

Let~$I_{k ,k}$ denote the~$k \times k$ identity matrix,
and
let $e^i_{k}\in \SA^k$   denote the $i$th canonical vector of size $k$, i.e.,
the~$i$th column of~$I_{k , k}$.
Let ${\mathcal L}^{k\times n}\subset \SA^{k\times n}$ denote the set of~$k\times n$ matrices
whose  columns are  all canonical vectors.

Using the semi-tensor product (\cite{cheng_book}) of matrices, denoted by~$\ltimes$,
the state-vector~$\begin{bmatrix} X_1(k) & \dots & X_n (k)\end{bmatrix}'$
 of a BCN is converted into
 a state-vector~$x(k) \in {\mathcal L}^{2^n}$.
 Basically,~$x(k)$ is the set of all the possible minterms of the~$X_i(k)$s,
 so~$x(k)$ is  a canonical vector for all~$k$.
 Similarly, the input vector~$\begin{bmatrix} U_1(k) & \dots & U_m (k)\end{bmatrix}'$
  is converted into
 a vector~$u(k) \in {\mathcal L}^{2^m}$.
 Since any Boolean  function   can be represented as a sum of minterms,
   the dynamics~\eqref{eq:BCN} can be represented in the bilinear form
            \begin{align} \label{eq:albac}
                        x(k+1)&=L \ltimes u(k)  \ltimes x(k).
            \end{align}
The matrix~$L \in {\mathcal L}^{2^n \times 2^{n+m}}$ is called the
\emph{transition matrix}  of the BCN.

 Algorithms for converting   a  BCN from the form~\eqref{eq:bcnio}
to its ASSR~\eqref{eq:albac}, and vice versa, may be found in~\cite{cheng_book}. Similarly,   the
BN~\eqref{eq:BNbinary}   may be represented  in  the  ASSR
            \begin{align} \label{eq:algebraic_BN}
                        x(k+1)&=L     x(k),
            \end{align}
where~$x(k)\in {\mathcal L}^{2^n}$ and~$L \in {\mathcal L}^{2^n \times 2^n}$.

The fact that a BN may be
represented in a linear form using the vector of minterms  has been known for a long time (see,
e.g.,~\cite{cull1971linear,cull1975control}), but the ASSR provides an explicit   algebraic form
that is particularly suitable for control-theoretic analysis.

 Given the ASSR~\eqref{eq:algebraic_BN} of a BN, we can associate        with it
a directed graph~$G =G(V,E)$, where $V=\{e^1_{2^n}, \dots, e^{2^n}_{2^n}\}$, and there
is a directed edge from vertex~$e^j_{2^n}$ to vertex~$e^i_{2^n}$ if and only
if~$[L]_{ij}=1$. In other words, there is a
 directed edge from vertex~$e^j_{2^n}$ to vertex~$e^i_{2^n}$
if and only if~$x(k)=e^j_{2^n}$ implies that~$x(k+1)=e^i_{2^n}$.

%%%%%%%%%%%%%%%%%%%%%%%%%%%%%%

We now briefly review some results from~\cite{hochma_sym}
derived by relating BCNs and
symbolic dynamics~(SD) (\cite{Lind_1995}).
SD has  evolved  from  analyzing general   dynamical
systems by discretizing the state-space into finitely many pieces,
  each labeled by a different symbol.
An orbit of the dynamical system is then transformed into
 a \emph{symbolic orbit} composed of
  the sequence of symbols corresponding
to the successive pieces  visited by the orbit.
The original evolution is transformed into a symbolic dynamics given by a  {shift
operator}~$\sigma$. The main object of study in SD is \emph{shift spaces}.

Given the BCN~\eqref{eq:bcnio},   define its \emph{set of state-trajectories of length~$j$}  by
  \begin{align*}
   {\mathcal A}_S^j&: = \{ X(0)  X(1)    \dots X(j-1)    :     \\& X(k+1)=f(X(k),U(k)), \;   U(k)\in \SA^m, \;  X(0)\in  \SA^n\},
  \end{align*}
 i.e.,  the   state trajectories  of length~$j$
 over  {all} possible controls and    initial
  conditions. Note that for a BN this becomes
  \begin{align*}
   \{ X(0)     \dots X(j-1)  :  \ X(k+1)=f(X(k) ),\;     X(0)\in \SA^n\}.
  \end{align*}
 %%%
The  {\em topological entropy} of a BCN   is  
%%%
\begin{align}\label{eq:tentBCN}
%%%
            h_S: =\lim_{j \to \infty} \frac{1}{j}\log |\mathcal{A}_S^j |   .
\end{align}
 %%%%%
%%%
In other words,~$h_S$
 is the asymptotic ``growth rate'' of
the number of state-sequences of a given length.
 A higher~$h_S$ corresponds to a ``richer'' control   in the sense
 that asymptotically more state-sequences  can be produced.

 \begin{example}\label{exa:minputBCN}
%%%%%%%%%%%%%%
Consider the     BCN:
\begin{align} \label{eq:minputBCN}
%%%
                X_1(k+1)&=U_1(k), \nonumber \\
                X_2(k+1)&=U_2(k), \nonumber \\
                 &\vdots \nonumber \\
                 X_m(k+1)&=U_m(k), \nonumber\\
                 X_{m+1}(k+1)&=f_1(X_1(k), \dots, X_n(k)) ,\\
                 X_{m+2}(k+1)&=f_2(X_1(k), \dots, X_n(k)) , \nonumber\\
                 &\vdots \nonumber\\
                X_n(k+1)&=f_{n-m}(X_1(k), \dots, X_n(k)).  \nonumber
%%%
\end{align}
It is straightforward to see that here~$|\mathcal{A}_S^j|=2^{n+(j-1)m}$,
so~\eqref{eq:tentBCN} yields
\be\label{eq:hsadded}
            h_S=\lim_{j\to \infty}\frac{1}{j} \left ((n+(j-1)m)\log 2 \right)=m\log 2.
\ee
Intuitively speaking, each
of the~$m$  control inputs in~\eqref{eq:minputBCN} contributes $\log2$ to the topological entropy.
%%%%%%%%%%%%
\end{example}

\cite{hochma_sym} have shown that in the ASSR,
 the set of state trajectories
                          of a BCN
                        is a  shift space (more precisely, a $1$-step shift space of finite type)  over the alphabet~$\{e_{2^n}^1,\dots ,
                        e_{2^n}^{2^n} \}$.
Combining this with
known results from SD   yields the following.
\begin{theorem}\label{thm:hsassr}(\cite{hochma_sym})
%%%%%%%%%%%%%%%%%%%%%%%%%%%%%%
Consider a BCN in the ASSR~\eqref{eq:albac}.
Let~$L_i:=L \ltimes e^i_{2^m}$, $i=1, \ldots , 2^m$,
 where~$L$ is the transition matrix of the BCN, and let
\be \label{def:M}
M := L_1 \vee L_2 \vee \ldots \vee L_{2^m}.
\ee
Then the topological entropy of the BCN  is
\be \label{eq:teproot}
                    h_S =\log \lambda_{M},
%%%%%%
\ee
where~$\lambda_{M}$ is the Perron root of the non-negative matrix~$M$.
\end{theorem}

\begin{remark}
%%%%%
Note that~$L_i \in \mathcal L^{2^n \times 2^{n }}$ and
thus~$ M \in \mathcal S^{2^n \times 2^n}$.
%%%%%
\end{remark}

\begin{example}\label{ex:entropybcn}
%%%%%%%%%%%%%%%%%%%%%%%
Consider the   BCN
defined by
\[
X(k+1)=U(k)  \vee  \bar{X}(k).
 \]
 The ASSR is given by~\eqref{eq:albac} with~$n=m=1$,
and $
            L=\begin{bmatrix}  1 & 1 & 0 &1 \\
             0 & 0 & 1 & 0   \end{bmatrix}
$. Thus, $
                    L_1 =L \ltimes e_2^1 = \begin{bmatrix} 1 & 1\\ 0  & 0 \end{bmatrix}$,
$                    L_2=L \ltimes e_2^2 = \begin{bmatrix} 0 & 1\cr 1 & 0 \end{bmatrix}$, and~$
M =L_1 \vee L_2=  \begin{bmatrix} 1 & 1\\ 1 & 0 \end{bmatrix}$. The
    eigenvalues of~$M$ are $  ({1 \pm \sqrt{5} })/2  $, so~\eqref{eq:teproot} yields
$h_S =\log ( ({1 +\sqrt{5} })/2)  $.
%%%%%%%%%%%
\end{example}

For easy reference, we recall the following result from the Perron-Frobenius
theory of non-negative matrices.
%%%%%%%%%%%%%%%%%%%%%%%%%%%%%%%%%%%%
\begin{theorem}\label{thm:min_max}~\cite[Ch.~8]{horn1}
 Suppose that~$A\in \mathbb{R}^{n\times n}_+$
 and let~$\lambda_{A}$ denote its Perron root. Then
  \begin{align} \label{eq:min_max2}
\min_{1\leq j\leq n}\sum_{i=1}^{n} A_{ij}  \leq \lambda_A \leq \max_{1\leq j\leq n} \sum_{i=1}^{n} A_{ij} .
\end{align}
 Furthermore,
 there exists~$w \in \R^n_+\setminus \{0\}$ such that~$A w =\lambda_{A} w $.
\end{theorem}

%%%%%%%%

%%%%%%%%%%%%%%%%%%%%%%%%%%%%%
\section{Main results}\label{sec:main}
%%%%%%%%%%%%%%%%%%%%%%%%%%%%%%%%%%%%%%%%%%%%%%%

Let~$\mathcal B_n^m$  denote
 the set of all BCNs with~$n$ state-variables and~$m$ inputs (with~$m\leq n$).
Let~$h_{max}$ be the maximum of the topological entropy over the BCNs in~$\mathcal B_n^m$.
Let~$  \mathcal {\overline B}_n^m \subset \mathcal B_n^m$
denote the subset of BCNs  with topological entropy equal to~$h_{max}$.
  A natural question is: what is the structure of the BCNs in~$   \mathcal {\overline B}_n^m$?

 Our first result shows in particular that the BCN in~\eqref{eq:minputBCN}
is in~$\mathcal {\overline B}_n^m$.
%%%
\begin{proposition} \label{prop:hm_max}
%%%
The maximal topological entropy of a~BCN in~$\mathcal B_n^m$ is $h_{max} = m\log 2$.
%%%
\end{proposition}

\noindent {\sl Proof.}
%%%%%%%%%%%%%%%%%%%%%%%%%%%%%%%%%%%%%%%%%%%%%%%%%
Fix a BCN in~$\mathcal B^m_n$, and consider its~ASSR.
Since~$M=\vee_{i=1}^{2^m}L_i$,
and every~$L_i $ has a single one entry in every column,
every column of~$M$ has no more than~$2^m$ one entries.
By~\eqref{eq:min_max2},~$\lambda_{M}\leq2^m$ so~$h_{S}\leq m\log 2$.
The BCN~\eqref{eq:minputBCN} attains  this bound
and this completes the proof.~$\square$

%%%%

Combining Theorems~\ref{thm:hsassr} and~\ref{thm:min_max} suggests that we can relate
the topological entropy of a BCN with the
maximum of the column (or row) sums of the matrix~$M$.
The next result shows that this is indeed so.
Let~$\alpha_{k,k}$ denote the~$k\times k$ matrix with all entries equal to~$\alpha$.
We use~$\alpha_k$ as a shorthand for~$\alpha_{k,1}$.

%%%%%%%%%%%%%%%%%%%%%%%%%%%%%%%%%%%%%%%%%%%%%%%%%%%%%%%%%%%%%%%%
\begin{proposition} \label{prop:rh_max}
%%%%%%%%%%%%%%%%%%%%%%%%%%%%%%%%%%%%%%%%%%%%%%%%%%%%%%%%%%%
 Consider a   BCN in the ASSR~\eqref{eq:albac}.
  Let
  \be\label{eq:defv}
  v:=\max_{1\leq j\leq 2^n}\sum_{i=1}^{ 2^n} M_{ij}  ,
  \ee
  where~$M$ is the matrix defined in~\eqref{def:M}.
   Then the following two conditions are equivalent.
  \begin{enumerate}[(a)]
\item \label{cond:a}  $h_S=\log v$. \label{cond:one}\\
\item \label{cond:b} There exist a permutation matrix $P\in \{0,1\}^{2^n\times 2^n}$
and~$r\geq v$ such that
%%%%%%%%%%%%%%%%%%%%%%%%%
\be \label{eq:rpred}
%%%%%%%%%%%%%%%%%%%%%%%%%%%%%%%%%%
P M P'=\begin{bmatrix}
B & C \\
0_{2^n-r, r} & D%
\end{bmatrix},
\ee
where~$B\in \mathcal S^{r\times r}$,     each column of~$B$ has
exactly~$v$ non zero elements, $D\in  \mathcal S ^{(2^n-r)\times
(2^n-r)}$, and~$C\in \mathcal S^{r\times (2^n-r)}$. \label{cond:two}
 \end{enumerate}
%%%%%%%%%%%%%%%%%%%%%%%%%%%%%%%%%%%%%%%%%%%%%%
\end{proposition}

{\sl Proof.}
%%%%%%%%%%%%%%%%%%%%%%%%%%%%%%%%%%%%%%%%%%
Assume that condition~\eqref{cond:two} holds.
Let~$w \in \R^r_+$ denote an eigenvector of~$B$ corresponding to its Perron root~$\lambda_B$.
Let~$\bar w:=\begin{bmatrix} w \\ 0_{2^{n}-r}  \end{bmatrix} $.
Then
\begin{align*}
                            P M P' \bar w &= \begin{bmatrix} B & C \\0_{2^n-r, r} & D \end{bmatrix} \bar w\\
                                         &=\lambda_B \bar w.
\end{align*}
This implies that~$P' \bar w $ is an eigenvector
of~$M$  corresponding to the eigenvalue~$\lambda_B$.
Since every column of~$B$ has exactly~$v$ one entries, Theorem~\ref{thm:min_max}
implies that~$\lambda_B=v$.
Combining this with~\eqref{eq:defv} and Theorem~\ref{thm:min_max}
implies that~$\lambda_M=v$, so~$h_S=\log \lambda_M=\log v$.
This shows that condition~\eqref{cond:two}
implies condition~\eqref{cond:one}.

To prove the converse implication, assume that~$h_S=\log v$. Then~$\lambda_M=v$.
By Theorem~\ref{thm:min_max}, there exists a vector~$w \in \R^{2^n}_+ \setminus \{0\}$
such that
$
M w= v  w   .
$
%%%%%
Let~$r \geq 1$ be the number of entries in~$w$ that are strictly positive, and
let~$P$ be a permutation matrix  such that
 \begin{align}\label{eq:deftildew}
 \tilde{w} & :=P w \nonumber\\
 &=\begin{bmatrix} \tilde{w}_1 &\tilde{w}_2& \dots& \tilde{w}_r& 0&\dots& 0\end{bmatrix}'
 \end{align}
 (note that if~$r=2^n$ then this vector includes no zeros).
 Then
 \be \label{eq:mwv}
                \tilde M \tilde w = v \tilde w,
 \ee
 where~$\tilde M:=P M P'$. Multiplying this on the left by~$1_{2^n}'$ yields
 \be \label{eq:sumtilde}
            \tilde s_1 \tilde{w}_1+\dots +\tilde s_r \tilde{w}_r  =    v (\tilde{w}_1+\dots +\tilde{w}_r ),
 \ee
 where~$\tilde s_i$ denotes the sum  of the elements in column~$i$ of~$\tilde M$.
 By~\eqref{eq:defv},~$\tilde s_i \leq v$ for all~$i$,
  so~\eqref{eq:sumtilde} implies that
 \be\label{eq:justadd}
                      \tilde s_i =v  \quad \text{for all }       i\in\{1,\dots,r \}.
 \ee
 Let
 $
                \tilde M= \begin{bmatrix} \tilde M_1 & \tilde M_2 \\
                \tilde M_3   & \tilde M_4      \end{bmatrix} ,
 $
 where~$\tilde M_1 \in \mathcal S^{r \times r}$.
%%%%
 Then~\eqref{eq:mwv} becomes
 \[
                \begin{bmatrix} \tilde M_1 & \tilde M_2 \\ \tilde M_3   & \tilde M_4      \end{bmatrix}
                \begin{bmatrix}\tilde w_1 & \dots & \tilde w_r &0 &\dots & 0 \end{bmatrix}'   =
                v\begin{bmatrix}\tilde w_1 & \dots & \tilde w_r &0 &\dots & 0 \end{bmatrix}'
                       .
 \]
%%%%%%%
Since the~$\tilde w_i$s are strictly positive, we conclude that~$\tilde M_3= 0_{2^n-r , r}$. Thus,~\eqref{eq:justadd} implies that
every column of~$\tilde M_1$ has exactly~$v$ one  entries,
so condition~\eqref{cond:two} holds.~$\square$

\begin{remark}\label{rem:intu}
%%%%%%%%%%
We can provide an intuitive explanation of\eqref{eq:rpred} as follows.
For a state~$a \in \{e_{2^n}^1, \dots, e_{2^n}^{2^n} \}$, let
\[
        R( a ):=\{ L \ltimes  e_{2^m}^1  \ltimes a ,\dots,  L \ltimes  e_{2^m}^{2^m}  \ltimes a
         \},
\]
i.e., the reachable set from~$a$ in one time step.
By the definition of~$M$,~$|R( e_{2^n}^j)|$
is equal to the number of one entries in column~$j$ of~$M$.
Thus,~$v$ is the maximal cardinality of the one time step
reachable sets. Proposition~\ref{prop:rh_max}
asserts that the topological entropy is equal to~$\log v$
if and only if there exists
a set~$Y$ containing~$r \geq v$  states
such that~$|R(a) |=v$ for all~$a \in Y$,
and any transition from a state in~$Y$ is to a state in~$Y$.
 \end{remark}

\begin{figure}[t]
 \centering
 \input{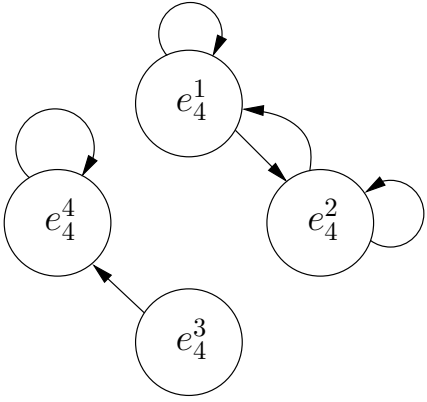}
\caption{State-space transition graph  of the BCN in Example~\ref{exa:rh_max}.}
\label{fig:state_space}
%\end{center}
\end{figure}

%%%%%%%%%%%%%%%%%%%%%%%
 \begin{example}\label{exa:rh_max}
 %%%%%%%%%%%%%%
 Consider the two-state, one-input BCN:
%%%%
\begin{align}   \label{eq:countr}
%%%%%
                    X_1(k+1)&= X_1(k),\\
                    X_2(k+1)&= [\bar{U} (k)\wedge  X_1(k)\wedge  \bar{X}_2(k)]\nonumber\\
                    & \vee [U (k)\wedge  X_1(k)\wedge  X_2(k)]. \nonumber
%%%%%
\end{align}
Fig.~\ref{fig:state_space} depicts the state-space transition graph of this~BCN, i.e.
a directed arrow from state~$a$ to  state~$b$ means that~$b$ belongs to the one time step reachable
set of~$a$.
It is easy to see from Fig.~\ref{fig:state_space} that~$v=2$
and that~$Y:=\{ e^1_4, e^2_4\}$
 satisfies the properties described in Remark~\ref{rem:intu}.
 By Proposition~\ref{prop:rh_max}, the topological entropy of~\eqref{eq:countr}
  is~$h_S=\log 2$.
%%%%%%%%%%%%%%%%%%%%%%%%%%%%%%%%%%%%%%%%%%%%%%%%%%%
%%%%%
\end{example}
%%%%%%%%%%%%%%%%%%%%%%%%%%%%%%%%%%%%%%%%%%%%%%%%%%%%%%%%%%%%%%%%%%
%%%%%%%%%%%%%%%%%%%%%%%%%%%%%%

\begin{example}
%%%
Consider again the
BCN in Example~\ref{exa:minputBCN}.
To analyze its topological entropy using
 Proposition~\ref{prop:rh_max} we
first derive an expression for
the matrix~$M$.

Let~$A \in \mathcal L^{2^{n-m} \times 2^n}$ denote
 the transition matrix in the~ASSR
with~$s:=n-m$ state-variables and~$m$ control inputs  BCN given by:
\begin{align*}
%%%
                  Y_{1}(k+1)&=f_{1}(W_1(k), \dots, W_m(k),Y_{ 1}(k),\dots, Y_{s}(k)),\\
                  &\vdots\\
                  Y_{s}(k+1)&=f_{s}(W_1(k), \dots, W_m(k),Y_{ 1}(k),\dots, Y_{s}(k)).
%%%
\end{align*}
%%%%
%%%%%%%%%%%%%%%%%%%%%%%%%%%%%%%%%%%%%%%
Pick~$ i \in \{ 1,\dots,2^m\}$. Consider the dynamics of~\eqref{eq:minputBCN} for~$u(k)=e^i_{2^m}$.
 By~\eqref{eq:minputBCN},~$x_1(k+1)\ltimes x_2(k+1) \ltimes \dots \ltimes x_m(k+1) = e^i_{2^m}$. Thus,
\begin{align*}
%%%
                 x(k+1)&=x_1(k+1) \ltimes \dots \ltimes x_m(k+1) \nonumber \\
                  &\;\;\;\;\;\ltimes x_{m+1}(k+1)\ltimes \dots \ltimes x_n(k+1)\nonumber\\
                 &=e^i_{2^m} \ltimes x_{m+1}(k+1)\ltimes \dots \ltimes x_n(k+1)\\
                 &=e^i_{2^m} \ltimes A  x(k) \nonumber\\
                 &=(e^i_{2^m}\otimes I_{2^{n-m},2^{n-m} })  A  x(k) \nonumber\\
                 &=\begin{bmatrix}
O_{(i-1) 2^{n-m},2^{n-m}} \\
I_{2^{n-m},2^{n-m} } \\
O_{(2^m-i)  2^{n-m},2^{n-m}}  \\
   \end{bmatrix}  A x(k). \nonumber
\end{align*}
%%%%
On the other-hand, for~$u(k)=e^i_{2^m}$,
 $x(k+1)=L_i x(k)$, so we conclude that
 $
            L_i=\begin{bmatrix}
O_{(i-1) 2^{n-m},2^n} \\
A \\
O_{(2^m-i)  2^{n-m},2^n}
   \end{bmatrix}  .
 $
Combining this with~\eqref{def:M} implies that
\begin{align} \label{eq:M_form}
M=\begin{bmatrix}  A\\
\vdots \\
A  \\
\end{bmatrix} \in \mathcal S^{2^n\times 2^n}.
\end{align}
Since every column of~$A$ is canonical vector,
 $M$ is a Boolean matrix and every column of~$M$ has exactly~$2^m$ ones.
Thus~$M$ has the form~\eqref{eq:rpred} with~$r=2^n$, $v=2^m$.
  Proposition~\ref{prop:rh_max}
  implies that~$h_S=m\log 2$,
and this agrees with~\eqref{eq:hsadded}.
%%%%%%%%%%%%%%
\end{example}
%%%%%%%%%%%%%%%%%%%%%%%

 One may perhaps expect that~\eqref{eq:minputBCN} is a ``canonical form''
  of a  BCN in~${\mathcal B^m_n}$, i.e.
    that for every BCN in this set there exists
     an invertible logical transformation
      of the state-variables taking it to the form~\eqref{eq:minputBCN}.
       However,  the next example
       shows that this is not so.
 \begin{example}\label{ex:contrexamle}
Consider  again the two--state, one-input~BCN in Example~\ref{exa:rh_max}.
Its ASSR is given by~$n=2$, $m=1$, and
\[
            L=\begin{bmatrix}  e_4^1 & e_4^2 & e_4^4 & e_4^4 & e_4^2 & e_4^1 & e_4^4 & e_4^4 \end{bmatrix}.
\]
 Thus, $ L_1 =L \ltimes e_2^1 = \begin{bmatrix} e_4^1 & e_4^2 & e_4^4 & e_4^4   \end{bmatrix}$,
$   L_2 =L \ltimes e_2^2 = \begin{bmatrix} e_4^2 & e_4^1 & e_4^4 & e_4^4   \end{bmatrix}$,
                     and
$M =L_1 \vee L_2=  \begin{bmatrix} e_4^1+ e_4^2 & e_4^1+ e_4^2 & e_4^4 & e_4^4   \end{bmatrix}$.
The   eigenvalues of~$M$ are~$\{ 2,1,0,0\}$,
 so $\lambda_M = 2 $ and~$h_S=\log 2$. Since~$2^m=2$,  Proposition~\ref{prop:hm_max} and
 Theorem~\ref{thm:hsassr} imply that~$h_S=h_{max}$.
 Since~$M$ has a unique zero row,~$P ' MP$ will
  also have a unique zero row, for any permutation matrix~$P$.
   Therefore~$P ' MP$ cannot have the form~\eqref{eq:M_form} for any permutation matrix~$P$.
%%%%%%%%%%%
%%%%%%%%%%%%
\end{example}
%%%%%%%%%%%%%%%%%%%%%%%%%%%%%%%%%%%%%%%%%%%%%%%%%%%

The next two results follow immediately from Proposition~\ref{prop:rh_max}.
%%%%%
 \begin{corollary} \label{coro:h_max}
A  BCN is in~$  \mathcal {\overline B}_n^m $  if and only if
condition~\eqref{cond:b} in Proposition~\ref{prop:rh_max}
holds and each column in the matrix~$B$
 has~$2^m$ non zero elements.
%%%%%%%%%%%%%%%%%%%%%%%%%%
\end{corollary}
%%%%%%%%%%%%%%%%%

\begin{corollary} \label{coro:M_max}
A  BCN is in~$\mathcal {\overline B}_n^n$ if and only if
\begin{align}\label{eq:fullm}
M=1_{2^n,2^n}.
\end{align}
%%%%%%%%%%%%%%%%%%%%%%%%%%%%%%%%%%%%%%%%
\end{corollary}

\begin{remark}
Recall that a BCN is called~$k$ \emph{fixed-time controllable}
if   for \emph{any}   $a,b\in  \{e_{2^n}^1, \dots,e_{2^n}^{2^n} \}  $
there exists a control that steers the BCN from~$x(0)=a$ to~$x(k)=b$
(see~\cite{dima_cont}).
%%%%%%%
Eq.~\eqref{eq:fullm} means that
 any state can be reached from any state in one time step.
 Thus, the BCN is $1$  fixed-time controllable.
\end{remark}

\subsection{Computational complexity}\label{subsec:comp}
%%%%%%%%%%%%%%%%%%%%%%%%%%%%%%%%%%%%%%%%%%%%%%%%%%%%%%%%%%%%%%%%%%%%%%%%
Consider the following problem.
\begin{problem}\label{prob:entropyproblem}
Given a BCN in~$ {\mathcal B^m_n}$
determine whether its topological entropy is~$h_S=h_{max}$.
\end{problem}
%%%%

\begin{proposition}\label{prop:hardness1}
 Problem~\ref{prob:entropyproblem} is NP-hard.
  \end{proposition}
%%%%%%%%%%%%%%%%%%%%%%%%%%

This implies that there does not exist  an algorithm with polynomial time complexity that
solves Problem~\ref{prob:entropyproblem}, unless~$P=NP$.

{\sl Proof  of Proposition~\ref{prop:hardness1}.}
%%%%%%%%%%%%%%%%%%%%%%%%%%%%%%%%%%%%%%%%%%%%%%%%%%%%%%
The proof  is based on a polynomial-time reduction of the famous
 SAT problem (see e.g.~\cite{garey1}) to Problem~\ref{prob:entropyproblem}.

 Consider a set of Boolean variables~$z_1,\dots,z_n$ taking values in~${\mathcal S}$.
 A \emph{formula} $g:{\mathcal S}^n \to {\mathcal S}$ is a rooted tree. The leaves include
 either a variable or its negation. Each internal node includes
 the operator~$\wedge$ or~$\vee$. The root of the tree then computes a
 formula in a natural way. The \emph{length} of the formula is the number of leaves in the tree.
 Formulas are often written as strings (e.g., $g(z_1,z_2)=(z_1 \wedge z_2)\vee \bar z_1)$),
  obtained by an inorder traversal of the rooted tree.

 A formula  is called \emph{satisfiable}
if there exists an assignment of its  variables
for which it attains the value~$1$.
 For example,~$g(z_1,z_2)=\bar{z}_1 \wedge  z_1 \wedge z_2$ is not satisfiable.
 %%%%%%%%%%%%%%%%%%%%%%%%%%%%%%%%%%%%%%%%%%%%%%%%%%%%%%%%%%%%%%%%%%%%%%%%%%%%%%%%%%%%%%%%%%%%%
 \begin{problem}\label{prob:sat}(SAT)
Given a Boolean formula~$g: {\mathcal S}^n \to {\mathcal S}$, determine whether it is satisfiable.
\end{problem}

Given a formula~$g:{\mathcal S}^n \to {\mathcal S}$, consider the~BCN in~${\mathcal B^n_n}$
defined by
\begin{align*}
                        X_1(k+1)&=U_1(k) \wedge (1-g(X_1(k),\dots,X_n(k)) ), \\
                        &\vdots\\
                        X_n(k+1)&=U_n(k)  \wedge (1-g(X_1(k),\dots,X_n(k))).
\end{align*}
It is clear that if~$g$ is not satisfiable then this  BCN is in~$  \mathcal {\overline B}_n^n $.
On the other-hand, if~$g$ is satisfiable then
there is at least one state
that is mapped to~$0_n$ for \emph{any} control.
This implies that in the ASSR,
at least one column of~$M$   is the vector~$e_{2^n}^{2^n}$.
Then Corollary~\ref{coro:M_max} implies that the~BCN is not in~$  \mathcal {\overline B}_n^n $.
 Summarizing, this provides a polynomial
 reduction from the~SAT problem to Problem~\ref{prob:entropyproblem}.
  Since~SAT is NP-complete even if the length of~$g$ is polynomial in~$n$,
  this completes the proof.~$\square$

 %%%%%%%%%%%%%%%%%%%%%%%%%%%%%%%%%%%%%%%%%%%%%%%%%%
%%%%%%%%%%%%%%%%%%%%%%%%%%%%%%%%%%%%%%%%%%%%%
\section{Conclusions} \label{sec:concl}
%%%%%%%%%%%%%%%%%%%%%%%%%%%%%%%%%%%%%%%5
BNs and BCNs are recently attracting considerable interest as computational models in systems
biology.

The topological entropy of a~BCN is a
 measure of how rich the control is.
 A   natural question is what is the structure of~BCNs
 that have the maximal possible topological entropy.
  In this paper,
  we derived  a necessary and sufficient condition for a~BCN to have this
   property, stated in terms of the~ASSR.

   Since the~ASSR of a~BCN with~$n$ state variables and~$m$
   inputs includes a matrix~$L \in {\mathcal L}^{2^n\times 2^{n+m}} $,
   verifying this conditions incurs an exponential time complexity.
   We also showed that the problem of determining whether a BCN has a
   maximal   topological entropy is NP-hard.
Thus, there does not exist an algorithm with polynomial time complexity that solves this problem, unless~$P=NP$.

Further research is needed in order to clarify the biophysical
 meaning  of
the topological entropy  in BCNs that model biological systems.
Another interesting  topic for further research is to characterize
 all
the possible values~$h$ such that there exists
a BCN in~${\mathcal B_n^m}$ with topological entropy~$h$.

\section*{Acknowledgements}
%%%%%%%%%%%%%%%%%%%%%
We thank the anonymous reviewers and the AE for their helpful comments.

 %%%%%%%%%%%%%%%%%%%%
 \bibliographystyle{natbib}
%%%%%%%%%%%%%%%%%%%%%%%%%%%
 \bibliography{Entrop_Dima}
%%%%%%%%%%
%\bibliographystyle{IEEEtranS} % this requires the file IEEEtranS.bst avialable on the web
 \end{document}